\documentclass[12pt,a4paper]{article}
\usepackage[cp1251]{inputenc} 
\usepackage{amsfonts}

\newcommand{\B}{$\hfill\Box$}
\newcommand{\al}{\alpha}

\newcommand{\de}{\delta}
\newcommand{\la}{\lambda}
\newcommand{\om}{\omega}

\newcommand{\vv}{\varphi}
\newcommand{\iy}{\infty}

\begin{document}

\begin{center}
{\large\bf
Recovering Differential Operators with Nonlocal Boundary Conditions}\\[0.2cm]
{\bf \dag Vjacheslav Anatoljevich Yurko and \ddag Chuan-Fu Yang} \\[0.2cm]
{\bf \dag Department of Mathematics, Saratov University,
Astrakhanskaya 83, Saratov 410012, Russia.
E-mail: yurkova@info.sgu.ru}\\[0.2cm]
{\bf \ddag Department of Applied Mathematics, Nanjing University of
Science and Technology, Nanjing, 210094, Jiangsu, China.
E-mail: chuanfuyang@njust.edu.cn}
\end{center}

\thispagestyle{empty}

{\bf Abstract.} Inverse spectral problems for Sturm-Liouville operators with nonlocal 
boundary conditions are studied. As the main spectral characteristics we introduce 
the so-called Weyl-type function and two spectra, which are generalizations of the well-known 
Weyl function and Borg's inverse problem for the classical Sturm-Liouville operator. Two 
uniqueness theorems of inverse problems from the Weyl-type function and two spectra
are presented and proved, respectively.

{\bf Key words:} differential operators; nonlocal boundary conditions;
inverse spectral problems

{\bf AMS Classification:} 34A55 34L05 47E05 \\

{\bf 1. Introduction. } Consider the differential equation
$$
-y''+q(x)y=\la y,\quad x\in (0,T),                                           \eqno(1)
$$
and the linear forms
$$
U_j(y):=\int_0^T y(t)d\sigma_j(t),\quad j=1,2.                               \eqno(2)
$$
Here $q(x)\in L(0,T)$ is a complex-valued function, $\sigma_j(t)$ are
complex-valued functions of bounded variations and are continuous from
the right for $t>0.$ There exist finite limits $H_j:=\sigma_j(+0)-
\sigma_j(0).$ Linear forms (2) can be written in the form
$$
U_j(y):=H_j y(0)+\int_0^T y(t)d\sigma_{j0}(t),\quad j=1,2,                   \eqno(3)
$$
where $\sigma_{j0}(t)$ in (3) are complex-valued functions of bounded
variations and are continuous from the right for $t\ge 0.$
For definiteness we assume that $H_1\ne 0.$

Boundary problems with nonlocal conditions are a part of fast developing
differential equations theory. Problems of this type arise in various fields of
physics, biology, biotechnology and etc. Nonlocal conditions come up when
value of the function on the boundary is connected to values inside the domain.
Theoretical investigation of problems with various type of nonlocal boundary
conditions is actual problem and recently it is paid much attention for them
in the literature. Originators of such problems were Samarskii and Bitsadze. 
They formulated and investigated nonlocal boundary problem for elliptic 
equation [1]. Afterwards the number of differential problems with nonlocal 
boundary conditions had increased.

Quite new area, related to problems of this type, deals with investigation
of the spectrum of differential equations with nonlocal conditions.
In this paper we study inverse spectral problems for Sturm-Liouville
operators with nonlocal boundary conditions, which are defined with
the help of linear forms (2) or (3). Classical inverse problems for
Eq. (1) with two-point separated boundary conditions have been studied
fairly completely in many works (see the monographs [2]-[5] and the
references therein). The theory of nonlocal inverse spectral problems
now is only at the beginning because of its complexity. Some aspects
of the inverse problem theory for different nonlocal operators can be
found in [6]-[15]. In this paper we prove uniqueness theorems for the
solution of the inverse spectral problems for Eq. (1) with nonlocal
boundary conditions. In Section 1 we suggest statements of the inverse
problems and formulate our main results (Theorems 1 and 2). Section 2
introduces important notions and properties of spectral characteristics.
The proofs of Theorems 1 and 2 are given in Section 3. For this purpose
we use the ideas of the method of spectral mappings [5]. In Section 4 we
present counterexamples related to the statements of the inverse problems
(see also [8]). Additional spectral data are introduced in Section 5.
In Section 6, as an example, we consider the inverse problem of recovering
the potential $q$ from the given three spectra.

\smallskip
Let $X_k(x,\la)$ and $Z_k(x,\la),$ $k=1,2,$ be the solutions of Eq. (1)
under the initial conditions
$$
X_1(0,\la)=X_2'(0,\la)=Z_1(T,\la)=Z_2'(T,\la)=1,
$$
$$
X_1'(0,\la)=X_2(0,\la)=Z_1'(T,\la)=Z_2(T,\la)=0.
$$
Consider the boundary value problem (BVP) $L_0$ for Eq. (1) with the conditions
$$
U_1(y)=U_2(y)=0.
$$
Denote $\om(\la):=\det[U_j(X_k)]_{j,k=1,2}$, and assume that $\om(\la)\not\equiv 0.$
The function $\om(\la)$ is entire in $\la$ of order $1/2,$ and its zeros
$\Xi=\{\xi_n\}_{n\geq 1}$ coincide with the eigenvalues of $L_0$. The function $\om(\la)$
is called the characteristic function for $L_0$.

Denote $V_j(y):=y^{(j-1)}(T),\; j=1,2.$ Consider the BVP $L_j$, $j=1,2,$ for
Eq. (1) with the conditions $U_j(y)=V_1(y)=0.$ The eigenvalues $\Lambda_j=
\{\la_{nj}\}_{n\ge 1}$ of the BVP $L_j$ coincide with the zeros of the
characteristic function $\Delta_j(\la):=\det[U_j(X_k), V_1(X_k)]_{k=1,2}$.

Let $\Phi(x,\la)$ be the solution of Eq. (1) under the conditions $U_1(\Phi)=1,$
$V_1(\Phi)=0.$ Denote $M(\la):=U_2(\Phi).$ The function $M(\la)$ is called the
Weyl-type function. It is known [4] that for Sturm-Liouville operators with
classical two-point separated boundary conditions, the specification of the
Weyl function uniquely determines the potential $q(x).$ In our case with
nonlocal boundary conditions, it is not true; the specification of the Weyl-type
function $M(\la)$ does not uniquely determine the potential (see counterexamples
in Section 4). In our case the inverse problem is formulated as follows.

\smallskip
{\bf Inverse problem 1.} Let $\Lambda_1\cap\Xi=\emptyset$ (condition $S$).
Given $M(\la)$ and $\om(\la),$ construct the potential $q(x).$

\smallskip
We note that the functions $\sigma_j(t)$ are known a priori, and only the
potential $q(x)$ has to be constructed.

Let us formulate a uniqueness theorem for Inverse problem 1. For this purpose,
together with $q$ we consider another potential $\tilde q,$ and we agree that
if a certain symbol $\al$ denotes an object related to $q,$ then $\tilde\al$
will denote an analogous object related to $\tilde q.$

\smallskip
{\bf Theorem 1. }{\it Let $\Lambda_1\cap\Xi=\emptyset$. If $M(\la)=\tilde M(\la)$
and $\om(\la)=\tilde\om(\la),$ then $q(x)=\tilde q(x)$ a.e. on $(0,T).$}

\smallskip
Thus, under condition $S,$ the specification $M(\la)$ and $\om(\la)$ uniquely
determines the potential. The proof of Theorem 1 see below in Section 3. We
note that if condition $S$ does not hold, then the specification $M(\la)$ and
$\om(\la)$ does not uniquely determine the potential (see counterexamples in
Section 4). In this case we have to specify an additional spectral information
(see Section 5).

\smallskip
Consider the BVP $L_{11}$ for Eq. (1) with the conditions $U_1(y)=V_2(y)=0.$
The eigenvalues $\Lambda_{11}:=\{\la_{n1}^1\}_{n\ge 1}$ of the BVP $L_{11}$
coincide with the zeros of the characteristic function $\Delta_{11}(\la):=
\det[U_1(X_k), V_2(X_k)]_{k=1,2}$. Clearly, $\{\lambda_{n1}\}_{n\geq 1}\cap\{\lambda_{n,1}^1\}_{n\geq 1}=\emptyset.$

\smallskip
{\bf Inverse problem 2. } Given $\{\la_{n1}, \la_{n1}^1\}_{n\ge 1}$,
construct $q(x).$

\smallskip
This inverse problem is a generalization of the well-known Borg's inverse
problem [16] for Sturm-Liouville operators with classical two-point
separated boundary conditions, and coincides with it when $U_1(y)=y(0).$ We note
that in Inverse problem 2 there are no restrictions on behavior of the spectra.

\smallskip
{\bf Theorem 2. }{\it If $\la_{n1}=\tilde\la_{n1}, \la_{n1}^1=\tilde\la_{n1}^1$,
$n\ge 1,$ then $q(x)=\tilde q(x)$ a.e. on $(0,T).$}

\smallskip
The proof of Theorem 2 see below in Section 3.\\

{\bf 2. Auxiliary propositions.} Let $\la=\rho^2,\; \tau :=\mbox{Im}\,\rho\ge 0.$
It is known (see, for example, [4]) that there exists a fundamental system of
solutions $\{Y_k(x,\rho)\}_{k=1,2}$ of Eq. (1) such that for $|\rho|\to\iy$:
$$
Y_1^{(\nu)}(x,\rho)\!=\!(i\rho)^{\nu}\exp(i\rho x)(1\!+\!O(\rho^{-1})),\;
Y_2^{(\nu)}(x,\rho)\!=\!(-i\rho)^{\nu}\exp(-i\rho x)(1\!+\!O(\rho^{-1})),              \eqno(4)
$$
$$
\det[Y_k^{(\nu-1)}(x,\rho)]_{k,\nu=1,2}=-2i\rho(1+O(\rho^{-1})).               \eqno(5)
$$

{\bf Lemma 1. }{\it Let $\{W_k(x,\la)\}_{k=1,2}$ be a fundamental system of
solutions of Eq. (1), and let $Q_j(y),\; j=1,2,$ be linear forms. Then}
$$
\det[Q_j(W_k)]_{k,j=1,2}=\det[Q_j(X_k)]_{k,j=1,2}
\det[W_k^{(\nu-1)}(x,\la)]_{k,\nu=1,2}\,.                                      \eqno(6)
$$

{\it Proof.} One has for $\nu=1,2$,
$$
W_\nu(x,\la)=\sum_{k=1}^2 A_{\nu k}(\la) X_k(x,\la),
$$
where the coefficients $A_{\nu k}(\la)$ do not depend on $x.$ This yields
$$
\det[Q_j(W_k)]_{k,j=1,2}=\det[Q_j(X_k)]_{k,j=1,2}\det[A_{\nu k}(\la)]_{k,j=1,2}\,,
$$
and
$$
\det[W_k^{(\nu-1)}(x,\la)]_{k,\nu=1,2}=\det[X_k^{(\nu-1)}(x,\la)]_{k,\nu=1,2}
\det[A_{\nu k}(\la)]_{k,j=1,2}\,.
$$
Since $\det[X_k^{(\nu-1)}(x,\la)]_{k,\nu=1,2}=1$, we arrive at (6).
\B

\medskip
It follows from (5)-(6) that
$$
\det[Q_j(Z_k)]_{k,j=1,2}=\det[Q_j(X_k)]_{k,j=1,2},                            \eqno(7)
$$
$$
\det[Q_j(Y_k)]_{k,j=1,2}=-2i\rho(1+O(\rho^{-1})) \det[Q_j(X_k)]_{k,j=1,2}\,.  \eqno(8)
$$

Consider the functions
$$
\vv(x,\la)=-\det[X_k(x,\la), U_1(X_k)]_{k=1,2},\;
\theta(x,\la)=\det[X_k(x,\la), U_2(X_k)]_{k=1,2},
$$
$$
\psi(x,\la)=\det[X_k(x,\la), V_1(X_k)]_{k=1,2}.
$$
Clearly,
$$
U_1(\vv)=0,\;U_2(\vv)=\om(\la),\;V_1(\vv)=\Delta_1(\la),\;
V_2(\vv)=\Delta_{11}(\la),
$$
$$
U_1(\theta)=\om(\la),\; U_2(\theta)=0,\; V_1(\theta)=-\Delta_2(\la),
$$
$$
U_j(\psi)=\Delta_j(\la),\; V_1(\psi)=0,\; V_2(\psi)=-1.
$$
Moreover, using (6)-(7), we calculate
$$
\det[\theta^{(\nu-1)}(x,\la), \vv^{(\nu-1)}(x,\la)]_{\nu=1,2}\!=\!\om(\la),\;
\det[\psi^{(\nu-1)}(x,\la), \vv^{(\nu-1)}(x,\la)]_{\nu=1,2}\!=\!\Delta_1(\la),    \eqno(9)
$$
$$
\Delta_1(\la)=-U_1(Z_2),\;\Delta_2(\la)=-U_2(Z_1),\;\Delta_{11}(\la)=U_1(Z_1). \eqno(10)
$$
Comparing boundary conditions on $\Phi, \psi, \vv$ and $\theta,$ we obtain
$$
\Phi(x,\la)=\frac{\psi(x,\la)}{\Delta_1(\la)},                                \eqno(11)
$$
$$
\Phi(x,\la)=\frac{1}{\om(\la)}\,\Big(\theta(x,\la)+\frac{\Delta_2(\la)}{\Delta_1(\la)}\vv(x,\la)\Big).     \eqno(12)
$$
Hence,
$$
M(\la):=U_2(\Phi)=\frac{\Delta_2(\la)}{\Delta_1(\la)},                                   \eqno(13)
$$
$$
\det[\Phi^{(\nu-1)}(x,\la), \vv^{(\nu-1)}(x,\la)]_{\nu=1,2}=1.                \eqno(14)
$$
Let $v_1(x,\la)$ and $v_2(x,\la)$ be the solutions of Eq. (1) under
the conditions
$$
v_1(T,\la)=v_2'(T,\la)=1,\quad v_1'(T,\la)=0,\quad U_1(v_2)=0.
$$
Obviously,
$$
v_1(x,\la)=Z_1(x,\la),\; v_2(x,\la)=Z_2(x,\la)+N(\la)Z_1(x,\la),
\; \det[v_k^{(\nu-1)}(x,\la)]_{k,\nu=1,2}=1,                                 \eqno(15)
$$
where
$$
N(\la)=\frac{\Delta_1(\la)}{\Delta_{11}(\la)}=-\frac{U_1(Z_2)}{U_1(Z_1)}\,.  \eqno(16)
$$

Denote
$$
U_1^a(y):=\int_0^a y(t)d\sigma_1(t),\quad a\in(0,T].
$$
Clearly, $U_1=U_1^T,$ and if $\sigma_1(t)\equiv C$ (constant) for $t\ge a,$
then $U_1=U_1^a.$

Let $\la_{n1}=\rho_{n}^2.$ For sufficiently small $\de>0,$ we denote
$$
\Pi_\de:=\{\rho:\; \mbox{arg}\,\rho\in[\de,\pi-\de]\},\ \
G_\de:=\{\rho:\; |\rho-\rho_n|\ge\de,\;\;\forall n\ge 1\}.
$$

{\bf Lemma 2. }{\it For $|\rho|\to\iy,\;\rho\in\Pi_\de$,
$$
\Phi^{(\nu)}(x,\la)=
\frac{(i\rho)^\nu}{H_1}\,\exp(i\rho x)(1+o(1)),\;x\in[0,T),                 \eqno(17)
$$
$$
v_1^{(\nu)}(x,\la)=\frac{(i\rho)^\nu}{2}
\,\exp(-i\rho(T-x))(1+O(\rho^{-1})),\; x\in[0,T),                           \eqno(18)
$$
$$
\Delta_1(\la)=-\frac{H_1}{2i\rho}\,\exp(-i\rho T)(1+o(1)),\;
\Delta_{11}(\la)=\frac{H_1}{2}\,\exp(-i\rho T)(1+o(1)).                     \eqno(19)
$$

Let $\sigma_1(t)\equiv C$ (constant) for $t\ge a$ (i.e. $U_1=U_1^a$). Then for
$|\rho|\to\iy,\;\rho\in\Pi_\de$,}
$$
\vv^{(\nu)}(x,\la)=\frac{H_1}{2}\,(-i\rho)^{\nu-1}
\exp(-i\rho x)(1+o(1)+O(\exp(i\rho(2x-a)))),\; x\in(0,T],                   \eqno(20)
$$
$$
v_2^{(\nu)}(x,\la)=(-i\rho)^{\nu-1}
\exp(i\rho(T-x))(1+o(1)+O(\exp(i\rho(2x-a)))),\; x\in[0,T).                 \eqno(21)
$$

{\it Proof.} One has
$$
\Phi(x,\la)=A_1(\la)Y_1(x,\rho)+A_2(\la)Y_2(x,\rho).                        \eqno(22)
$$
Since $U_1(\Phi)=1,\; V_1(\Phi)=0,$ it follows from (22) that
$$
A_1(\la)U_1(Y_1)+A_2(\la)U_1(Y_2)=1,\; A_1(\la)V_1(Y_1)+A_2(\la)V_1(Y_2)=0. \eqno(23)
$$
By virtue of (4), we have for $|\rho|\to\iy,\; \rho\in\Pi_\de$:
$$
U_1(Y_1)=H_1(1+o(1)),\; U_1(Y_2)=O(\exp(-i\rho T)),                         \eqno(24)
$$
$$
V_1(Y_1)=\exp(i\rho T)(1+O(\rho^{-1})),
\; V_1(Y_2)=\exp(-i\rho T)(1+O(\rho^{-1})).                                 \eqno(25)
$$
Solving linear algebraic system (23) and using (24)-(25), we calculate
$$
A_1(\rho)=H_1^{-1}(1+o(1)),\; A_2(\rho)=O(\exp(2i\rho T)).
$$
Substituting these relations into (22), we arrive at (17).
Formulas (18)-(21) are proved similarly.
\B

\smallskip
By the well-known method (see, for example, [4]) one can also obtain the
following estimates for $x\in(0,T),\; \tau\ge 0:$
$$
v_1^{(\nu)}(x,\la)=O(\rho^{\nu}\exp(-i\rho(T-x))),                         \eqno(26)
$$
$$
\Phi^{(\nu)}(x,\la)=O(\rho^{\nu}\exp(i\rho x)),\quad \rho\in G_\de.        \eqno(27)
$$
Moreover, if $\sigma_1(t)\equiv C$ (constant) for $t\ge a$ (i.e. $U_1=U_1^a$),
then for $x\ge a/2,\; \tau\ge 0$:
$$
\vv^{(\nu)}(x,\la)=O(\rho^{\nu-1}\exp(-i\rho x)),                          \eqno(28)
$$
$$
v_2^{(\nu)}(x,\la)=O(\rho^{\nu-1}\exp(i\rho(T-x))),\quad \rho\in G_\de.    \eqno(29)
$$

\bigskip
{\bf 3. Proofs of Theorems 1-2.} Firstly we prove Theorem 2. Let
$\la_{n1}=\tilde\la_{n1}, \la_{n1}^1=\tilde\la_{n1}^1$, $n\ge 1.$
The characteristic function $\Delta_1(\la)$ of the BVP $L_1$ is entire
in $\la$ of order $1/2.$ Therefore, by Hadamard's factorization theorem,
$\Delta_1(\la)$ is uniquely determined up to a multiplicative constant
by its zeros, i.e. $\Delta_1(\la)/\tilde\Delta_1(\la)\equiv C$ (constant).
Taking (19) into account, we calculate $C=1,$ and consequently,
$\Delta_1(\la)\equiv\tilde\Delta_1(\la).$ Analogously, we get
$\Delta_{11}(\la)\equiv\tilde\Delta_{11}(\la).$ By virtue of (16),
this yields
$$
N(\la)\equiv\tilde N(\la).                                               \eqno(30)
$$
 Consider the functions
$$
P_1(x,\la)\!\!=\!\!v_1(x,\la)\tilde v'_2(x,\la)\!-\!\tilde v'_1(x,\la)v_2(x,\la),
\;P_2(x,\la)\!\!=\!\!v_2(x,\la)\tilde v_1(x,\la)\!-\!\tilde v_2(x,\la)v_1(x,\la).    \eqno(31)
$$
In view of (15) and (30), one gets
$$
P_1(x,\la)=(Z_1(x,\la)\tilde Z'_2(x,\la)-\tilde Z'_1(x,\la)Z_2(x,\la))
+(\tilde N(\la)-N(\la))Z_1(x,\la)\tilde Z'_1(x,\la)
$$
$$
=Z_1(x,\la)\tilde Z'_2(x,\la)-\tilde Z'_1(x,\la)Z_2(x,\la),
$$
$$
P_2(x,\la)=Z_2(x,\la)\tilde Z_1(x,\la)-\tilde Z_2(x,\la)Z_1(x,\la)
+(N(\la)-\tilde N(\la))Z_1(x,\la)\tilde Z_1(x,\la)
$$
$$
=Z_2(x,\la)\tilde Z_1(x,\la)-\tilde Z_2(x,\la)Z_1(x,\la).
$$
Thus, for each fixed $x,$ the functions $P_k(x,\la),\; k=1,2,$
are entire in $\la.$ On the other hand, taking (18) and (21) into
account we calculate for each fixed $x\ge T/2$ and $k=1,2:$
$$
P_k(x,\la)-\de_{1k}=o(1),\; |\rho|\to\iy,\;\rho\in\Pi_\de,
$$
where $\de_{1k}$ is the Kronecker symbol. Moreover, in view of
(26) and (29), we get for $k=1,2$
$$
P_k(x,\la)=O(1),\; |\rho|\to\iy,\;\rho\in \Pi_\de.
$$
Using the maximum modulus principle and Liouville's theorem for
entire functions, we conclude that
$$
P_1(x,\la)\equiv 1,\quad P_2(x,\la)\equiv 0,\quad x\ge T/2.
$$
Together with (31) this yields
$$
v_k(x,\la)=\tilde v_k(x,\la),\;
Z_k(x,\la)=\tilde Z_k(x,\la),\; q(x)=\tilde q(x),\; x\ge T/2.           \eqno(32)
$$

Let us now consider the BVPs $L_1^a$ and $L_{11}^a$ for Eq. (1) on
the interval $(0,T)$ with the conditions $U_1^a(y)=V_1(y)=0$ and
$U_1^a(y)=V_2(y)=0,$ respectively. Then, according to (10), the
functions $\Delta_1^a(\la):=-U_1^a(Z_2)$ and $\Delta_{11}^a(\la)
:=U_1^a(Z_1)$ are the characteristic functions of $L_1^a$ and
$L_{11}^a,$ respectively. One has
$$
U_1^{a/2}(Z_k)=U_1^{a}(Z_k)-\int_{a/2}^a Z_k(t,\la)d\sigma_1(t),
\quad k=1,2,
$$
hence
$$
\Delta_{1}^{a/2}(\la)=
\Delta_{1}^a(\la)+\int_{a/2}^a Z_2(t,\la)d\sigma_1(t),\;
\Delta_{11}^{a/2}(\la)=
\Delta_{11}^a(\la)-\int_{a/2}^a Z_1(t,\la)d\sigma_1(t).                 \eqno(33)
$$
Let us use (33) for $a=T.$ Since $\Delta_{1}^T(\la)=\Delta_{1}(\la),\;
\Delta_{11}^T(\la)=\Delta_{11}(\la),$ it follows from (32)-(33) that
$$
\Delta_{1}^{T/2}(\la)=\tilde\Delta_{1}^{T/2}(\la),\quad
\Delta_{11}^{T/2}(\la)=\tilde\Delta_{11}^{T/2}(\la).
$$
Repeating preceding arguments subsequently for $a=T/2, T/4, T/8,\ldots,$
we conclude that $q(x)=\tilde q(x)$ a.e. on $(0,T).$ Theorem 2 is proved.
\B

\smallskip
Now we will prove Theorem 1. Let $\Lambda_1\cap\Xi=\emptyset,$ and let
$M(\la)=\tilde M(\la),$ $\om(\la)=\tilde\om(\la).$ Consider the functions
$$
R_1(x,\la)\!=\!\Phi(x,\la)\tilde\vv'(x,\la)\!-\!\tilde\Phi'(x,\la)\vv(x,\la),
\;R_2(x,\la)\!=\!\vv(x,\la)\tilde\Phi(x,\la)\!-\!\tilde\vv(x,\la)\Phi(x,\la).   \eqno(34)
$$
It follows from (11) and (34) that
$$
R_1(x,\la)=\frac{1}{\Delta_1(\la)}
\Big(\psi(x,\la)\tilde\vv'(x,\la)-\tilde\psi'(x,\la)\vv(x,\la)\Big),
$$
$$
R_2(x,\la)=\frac{1}{\Delta_1(\la)}
\Big(\vv(x,\la)\tilde\psi(x,\la)-\tilde\vv(x,\la)\psi(x,\la)\Big).
$$
This yields that for each fixed $x,$ the functions $R_k(x,\la)$ are
meromorphic in $\la$ with possible poles only at $\la=\la_{n1}$.
On the other hand, taking (12) into account we calculate
$$
R_1(x,\la)=\frac{1}{\om(\la)}
\Big(\theta(x,\la)\tilde\vv'(x,\la)-\tilde\theta'(x,\la)\vv(x,\la)\Big), \eqno(35)
$$
$$
R_2(x,\la)=\frac{1}{\om(\la)}
\Big(\vv(x,\la)\tilde\theta(x,\la)-\tilde\vv(x,\la)\theta(x,\la)\Big).   \eqno(36)
$$
Hence the functions $R_k(x,\la)$ are regular at $\la=\la_{n1}$. Thus,
for each fixed $x,$ the functions $R_k(x,\la)$ are entire in $\la.$
Furthermore, by virtue of (17) and (20), we obtain for $x\ge T/2:$
$$
R_k(x,\la)-\de_{1k}=o(1),\quad |\rho|\to\iy,\; \rho\in\Pi_\de.
$$
Moreover, using (27)-(28), we get for $x\geq T/2:$
$$
R_k(x,\la)=O(1),\quad |\rho|\to\iy,\; \rho\in G_\de.
$$
Therefore, $R_1(x,\la)\equiv 1,\; R_2(x,\la)\equiv 0.$
Together with (14) and (34), this yields
$$
\vv(x,\la)=\tilde\vv(x,\la),\; \psi(x,\la)=\tilde\psi(x,\la),\;
q(x)=\tilde q(x),\; x\ge T/2.
$$
In particular, we obtain
$$
Z_k(x,\la)=\tilde Z_k(x,\la),\quad k=1,2,\quad x\ge T/2.
$$
Since
$$
\vv(x,\la)=U_1(Z_1)Z_2(x,\la)-U_1(Z_2)Z_1(x,\la),
$$
it follows that
$$
\Delta_{1}(\la)=\tilde\Delta_{1}(\la),\;\Delta_{11}(\la)=\tilde\Delta_{11}(\la).
$$
Using Theorem 2, we conclude that $q(x)=\tilde q(x)$ a.e. on $(0,T).$
Theorem 1 is proved.
\B

\bigskip
{\bf 4. Counterexamples. }
1) Let $T=\pi,\;U_1(y)=y(0),\;U_2(y)=y(\pi/2),$ $q(x)=q(x+\pi/2),\;x\in(0,\pi/2),$
and $q(x)\not\equiv q(\pi-x).$ Take $\tilde q(x)=q(\pi-x),\; x\in (0,\pi).$ We see that the BVP
$\tilde{L}_1$, for
Eq. (1) with $\tilde{q}(x)=q(\pi-x)$ and the conditions $U_1(y)=V_1(y)=0$;
the BVP
$\tilde{L}_2$, for
Eq. (1) with $\tilde{q}(x)=q(\pi-x)$ and the conditions $U_2(y)=V_1(y)=0$;
and the BVP
$\tilde{L}_0$, for
Eq. (1) with $\tilde{q}(x)=q(\pi/2-x)$ and the conditions $U_1(y)=U_2(y)=0$.

Then
$$
\Delta_1(\la)=\tilde\Delta_1(\la),\; \Delta_2(\la)=\tilde\Delta_2(\la),\;
\om(\la)=\tilde\om(\la),
$$
and, in view of (13), $M(\la)=\tilde M(\la).$ Condition $S$ does not hold.
This means, that the specification of $M(\la)$ and $\om(\la)$ does not
uniquely determine the potential $q.$

\smallskip
2) Let $T=\pi,\;U_1(y)=y(0),\;U_2(y)=y(\pi-\al),$ where $\al\in(0,\pi/2).$ Then
$$
\Delta_1(\la)=X_2(\pi,\la),\; \om(\la)=X_2(\pi-\al,\la),
$$
and
$$
\Delta_2(\la)=X_2(\pi,\la)X_1(\pi-\al,\la)-X_2(\pi-\al,\la)X_1(\pi,\la).
$$
Obviously, if $\Delta_1(\la^*)\Delta_2(\la^*)\om(\la^*)=0$ for a certain $\la^*$,
then either $\Delta_1(\la^*)=\Delta_2(\la^*)=\om(\la^*)=0$ (i.e. $\la^*$ is an
eigenvalue for all boundary value problems $L_0, L_1, L_2$), or $\la^*$ is an
eigenvalue for only one problem from $L_0, L_1, L_2$. In other words, it is
impossible that $\la^*$ is an eigenvalue for only two problems from $L_0, L_1, L_2$.

Let $q(x)\not\equiv q(\pi-x),$ and let $q(x)\equiv 0$ for $x\in[0,\al_0]\cup
[\pi-\al_0,\pi],$ where $\al_0\in(0,\pi/2).$ If $\al<\al_0$, then $\la_{n2}=
(\pi n/\al)^2,\; n\ge 1.$ Choose a sufficiently small $\al<\al_0$ such that
$\Lambda_1\cap\Lambda_2=\emptyset.$ Clearly, such choice is possible. Then
$\Lambda_1\cap\Xi=\emptyset,$ i.e. condition $S$ holds. Take
$\tilde q(x):=q(\pi-x).$ Then $\Delta_1(\la)=\tilde\Delta_1(\la),$
$\Delta_2(\la)=\tilde\Delta_2(\la),$ and consequently, $M(\la)=\tilde M(\la).$
Thus, condition $S$ holds, but the specification of $M(\la)$ does not uniquely
determine the potential $q.$

\bigskip
{\bf 5. Additional spectral data.}
If condition $S$ does not hold, then the specification of $M(\la)$ and
$\om(\la)$ does not uniquely determine the potential $q.$ We introduce
additional spectral data. For simplicity, we confine ourselves to the
case when zeros of $\om(\la)$ are simple. By virtue of (9),
$$
\det[\theta^{(\nu-1)}(x,\la), \vv^{(\nu-1)}(x,\la)]_{\nu=1,2}=\om(\la).
$$
Then the functions $\vv(x,\xi_n)$ and $\theta(x,\xi_n)$ are linearly dependent,
i.e. there exist numbers $A_n$ and $B_n$ ($|A_n|+|B_n|>0$) such that
$A_n\vv(x,\xi_n)=B_n\theta(x,\xi_n).$ Consider the sequence $D=\{d_n\}_{n\ge 1}$,
where $d_n:=B_n/A_n$ ($d_n:=\iy$ if $A_n=0$). The inverse problem is formulated
as follows.

\smallskip
{\bf Inverse problem 3. } Given $M(\la), \om(\la)$ and $D,$ construct $q(x).$

\smallskip
We note that if condition $S$ holds (i.e. $\Lambda_1\cap\Xi=\emptyset$), then
by virtue of (11)-(12), $d_n=-M^{-1}(\xi_n).$ At this case $M(\la)=\tilde{M}(\la)$
implies that $D=\tilde{D}$, and we arrive at Inverse problem 1.

\smallskip
{\bf Theorem 3. }{\it If $M(\la)=\tilde M(\la),$ $\om(\la)=\tilde\om(\la)$ and
$D=\tilde D,$ then $q(x)=\tilde q(x)$ a.e. on $(0,T).$}

\smallskip
{\it Proof. } Fix $n\ge 1,$ and consider the functions $R_k(x,\la),\; k=1,2,$
in a neighborhood of the point $\la=\xi_n$. If $d_n=\tilde d_n\ne\iy,$ then,
in view of (35)-(36), $R_k(x,\la)$ are regular at $\la=\xi_n$. If
$d_n=\tilde d_n=\iy,$ then $\theta(x,\xi_n)=\tilde\theta(x,\xi_n)=0.$
By virtue of (35)-(36), this yields that $R_k(x,\la)$ are regular at
$\la=\xi_n$. Thus, the functions $R_k(x,\la),\; k=1,2$ are entire in $\la.$
Repeating the arguments from the proof of Theorem 1, we conclude that
$q(x)=\tilde q(x)$ a.e. on $(0,T).$
\B

\bigskip
{\bf 6. Example (Inverse problem from three spectra).}
Fix $a\in(0,T).$ Consider Inverse problem 1 in the particular case, when
$U_1(y)=y(0),\; U_2(y)=y(a).$ Then the boundary value problems $L_0, L_1, L_2$
take the form
$$
L'_0:\qquad y(0)=y(a)=0,
$$
$$
L'_1:\qquad y(0)=y(T)=0,
$$
$$
L'_2:\qquad y(a)=y(T)=0.
$$

Denote by $\Lambda'_j=\{\la'_{nj}\}$ the spectrum of $L'_j,$ and assume that
$\Lambda'_0\cap\Lambda'_1=\emptyset$ (condition $S'$).

\smallskip
{\bf Inverse problem 4. } Given three spectra $\Lambda'_0, \Lambda'_1$ and
$\Lambda'_2$, construct $q(x).$

\smallskip
The following theorem is a consequence of Theorem 1.

\smallskip
{\bf Theorem 2. }{\it Let condition $S'$ hold. If $\Lambda'_j=\tilde\Lambda'_j$,
$j=0,1,2,$ then $q(x)=\tilde q(x)$ a.e. on $(0,T).$}

\smallskip
We note that Inverse problem 4 was studied by many authors
(see, for example, [17]-[18]).

\medskip
{\bf Acknowledgment.}  This research work of Yurko was supported in part by Grant
13-01-00134 of Russian Foundation for Basic Research. Yang was supported in
part by the National Natural Science Foundation of China (11171152) and
Natural Science Foundation of Jiangsu Province of China (BK 2014021904).

\begin{center}
{\bf REFERENCES}
\end{center}
\begin{enumerate}

\item[{[1]}] Bitsadze A.V. and Samarskii A.A., Some elementary generalizations of linear
elliptic boundary value problems. Dokl. Akad. Nauk SSSR 185 (1969), 739-740.

\item[{[2]}] Marchenko V.A., Sturm-Liouville operators and their applications.
     Naukova Dumka,  Kiev, 1977;  English  transl., Birkh\"auser, 1986.
\item[{[3]}] Levitan B.M., Inverse Sturm-Liouville problems. Nauka,
     Moscow, 1984; English transl., VNU Sci. Press, Utrecht, 1987.
\item[{[4]}] Freiling G. and Yurko V.A., Inverse Sturm-Liouville
     Problems and their Applications. NOVA Science Publishers, New York, 2001.
\item[{[5]}] Yurko V.A. Method of Spectral Mappings in the Inverse Problem Theory.
     Inverse and Ill-posed Problems Series. VSP, Utrecht, 2002.
\item[{[6]}] Yurko V.A. An inverse problem for integral operators. Matem. Zametki, 37,
     no.5 (1985), 690-701 (Russian); English transl. in  Mathematical Notes, 37,
		 no.5-6 (1985), 378-385.
\item[{[7]}] Yurko V.A. An inverse problem for integro-differential operators.
     Matem. Zametki, 50, no.5 (1991), 134-146 (Russian); English transl. in
		 Mathematical Notes, 50, no.5-6 (1991), 1188-1197.
\item[{[8]}] Kravchenko K.V. On differential operators with nonlocal boundary
     conditions. Differ. Uravn. 36, no.4 (2000), 464-469; English transl. in
		 Differ. Equations 36, no.4 (2000), 517-523.
\item[{[9]}] Buterin S.A. The inverse problem of recovering the Volterra convolution
     operator from the incomplete spectrum of its rank-one perturbation.
		 Inverse Problems 22 (2006), 2223-2236.
\item[{[10]}] Buterin S.A. On an inverse spectral problem for a convolution
     integro-differential operator. Results in Mathematics 50 (2007), no. 3/4, 173-181.
\item[{[11]}] Hryniv R., Nizhnik L.P., and Albeverio S. Inverse spectral problems for
     nonlocal Sturm-Liouville operators. Inverse Problems 23 (2007), 523-535.
\item[{[12]}] Nizhnik L.P. Inverse nonlocal Sturm-Liouville problem. Inverse Problems
     26, no.12 (2010), 125006, 9pp.
\item[{[13]}] Kuryshova Y.V. and Chung-Tsun Shieh. Inverse nodal problem for
     integro-differential operators. Journal of Inverse and Ill-Posed Proplems
		 18, no.4 (2010), 357-369.
\item[{[14]}] Freiling G. and Yurko V.A. Inverse problems for differential operators
     with a constant delay. Applied Mathematical Letters 25, no.11 (2012), 1999-2004.
\item[{[15]}] Yang C. F. Trace and inverse problem of a discontinuous Sturm-Liouville
     operator with retarded argument. J. Math. Anal. Appl. 395 (2012), no.1, 30-41.
\item[{[16]}] Borg G. Eine Umkehrung der Sturm-Liouvilleschen Eigenwertaufgabe.
     Acta Math. 78 (1946), 1-96.
\item[{[17]}] Gesztesy F., and Simon B. On the determination of a potential from three
     spectra. Amer. Math. Soc. Transl. Ser.2, 189, 85-92, Amer. Math. Soc.,
		 Providence, RI, 1999.
\item[{[18]}] Pivovarchik V. An inverse Sturm-Liouville problem by three spectra.
     IEOT 34 (1999), no.2, 234-243.
\end{enumerate}

\end{document}